\documentclass[10pt,preprint]{elsarticle}
\usepackage{amssymb,amsfonts,amsmath,mathrsfs,mathtools}
\usepackage[margin=1in,papersize={8.5in,11in}]{geometry}
\usepackage{xcolor}
\usepackage{multirow}
\usepackage{braket}
\usepackage{cleveref}
\usepackage{ulem}
\def\TEBD2{$\text{TEBD}^2$}
\def\bigO{\mathcal{O}}

\title{Computing excited states with isometric tensor networks in two-dimensions}

\begin{document}
\begin{frontmatter}

\author[lbnl]{Alec Dektor\corref{correspondingAuthor}}
\ead{adektor@lbl.gov}
\author[caltech]{Runze Chi}
\author[lbnl]{Roel Van Beeumen}
\author[lbnl]{Chao Yang}

\address[lbnl]{Applied Mathematics \& Computational Research Division, Lawrence Berkeley National Laboratory, Berkeley, California 94720, USA.}

\address[caltech]{Division of Chemistry and Chemical Engineering, California Institute of Technology, Pasadena, California 91125, USA}

\cortext[correspondingAuthor]{Corresponding author}

\journal{Arxiv}
\date{}

\begin{abstract}
We present a new subspace iteration method for computing low-lying eigenpairs (excited states) of high-dimensional quantum many-body Hamiltonians with nearest neighbor interactions on two-dimensional lattices. The method is based on a new block isometric projected entangled pair state (block-isoPEPS) ansatz that generalizes the block matrix product state (MPS) framework, widely used for Hamiltonians defined on one-dimensional chains, to two-dimensions. The proposed block-isoPEPS ansatz offers several attractive features for PEPS-based algorithms, including exact block orthogonalization, controlled local truncation via singular value decompositions, and efficient evaluation of observables. We demonstrate the proposed inexact subspace iteration for block-isoPEPS by computing excitations of the two-dimensional transverse-field Ising and Heisenberg models and compare our results with existing PEPS methods. Our results demonstrate that block isometric tensor networks provide a scalable framework for studying excitations in quantum many-body systems beyond one dimension. 

\end{abstract}


\end{frontmatter}

\section{Introduction} \label{sec:intro}

We consider the problem of computing a few of the algebraically smallest eigenpairs, i.e., ground and excited states, of Hamiltonians with pairwise interactions 
\begin{equation} \label{eq:Ham}
    H = \sum_{(j_1,j_2) \in S} 
    H_{j_1} H_{j_2}, \qquad 
    H_j = I^{\otimes(j-1)} \otimes M_j \otimes I^{\otimes(L-j)}
\end{equation} 
where $S \subset \{1,\ldots,L\}\times \{1,\ldots, L\}$ specifies the interacting pairs on $L$ local $d$-dimensional vector spaces (referred to as sites), each $M_j \in \mathbb{C}^{d\times d}$ is a Hermitian matrix, and $I$ denotes the $d \times d$ identity matrix. Hamiltonians of the form \eqref{eq:Ham} appear in model quantum many-body systems such as spin chains and lattices with nearest-neighbor interactions. The number of entries in the matrix  \eqref{eq:Ham} grows exponentially fast with the parameter $L$ making classical methods for finding eigenpairs too expensive for problems of interest. 

Tensor networks are a powerful framework for computing eigenvectors of \eqref{eq:Ham} using a reasonable number of degrees of freedom determined by an auxiliary parameter called the bond-dimension. For instance, when \eqref{eq:Ham} represents a Hamiltonian describing a quantum spin chain with nearest-neighbor interactions, i.e., $S$ contains only pairs of the form $(j,j+1)$, and has separated eigenvalues, the eigenvectors corresponding to the algebraically smallest eigenvalues admit accurate matrix product state (MPS) approximations with small bond-dimension \cite{Hastings2007}. 
Several MPS algorithms have been developed to compute eigenvectors \cite{White1993, Schollwock2011, zauner2018variational, motruk2016density, chepiga2017excitation, stoudenmire2012, zauner2018, Dektor2025, MPS_Lanczos} of \eqref{eq:Ham}. These algorithms enable practitioners to routinely compute accurate eigenvector approximations of (quasi) one-dimensional quantum many-body Hamiltonians \eqref{eq:Ham} with $L>100$. 

However, due to the one-dimensional connectivity of the MPS tensor network, these methods are often not suitable for quantum many-body systems with nearest-neighbor interactions in two-dimensional geometries such as rectangular lattices. In this case the set $S$ of interacting pairs contains $(j_1,j_2)$ with $j_1 \ll j_2$ which induces long-range interactions in the MPS, leading to prohibitively large bond-dimensions in MPS eigenvector approximations. These limitations have motivated the development of algorithms based on tensor networks with additional connectivity designed to efficiently capture the entanglement structure of quantum many-body systems defined on two- and higher-dimensional geometries. 

Projected entangled pair states (PEPS) are a class of tensor networks with connectivity designed to efficiently capture entanglement in two-dimensional geometries \cite{Orus2014}. Their connectivity makes PEPS significantly more demanding than MPS in terms of algorithmic design and computational resources. For example, exact contraction (required for computing norms, inner products, and observables) of PEPS networks has computational cost scaling exponentially with the 2D system size, whereas the corresponding contraction for MPS scales linearly with 1D system size (and polynomially with the bond-dimension). To mitigate the exponential scaling of PEPS contraction, several approximate schemes have been developed that trade accuracy for reduced cost \cite{Lubasch2014}. While these methods avoid exponential scaling, they remain computationally expensive, often scaling as high powers of the PEPS bond-dimension $\chi$, e.g., $\bigO(\chi^{10})$. 

An alternative strategy for obtaining computationally efficient PEPS algorithms is to impose isometric constraints on the tensors in the PEPS network, so that contracting large regions of the network is trivial. Imposing sufficient isometric constraints yields an isometric PEPS (isoPEPS) \cite{zaletel2020isometric} with an orthogonality center similar to canonical forms of MPS (which are central to many MPS algorithms). 
However unlike MPS, the orthogonality center cannot be moved around the PEPS without either increasing bond-dimension or introducing error. As such, 
the set of isoPEPS with fixed bond dimension $\chi$ is a strict subset of all PEPS with bond-dimension $\chi$. While the variational power of isoPEPS is still not fully understood, recent work has identified classes of physical states isoPEPS can represent exactly or efficiently \cite{Soejima2020} and has established complexity-theoretic limits on their contraction and sampling \cite{Malz2025}. The computational cost of many isoPEPS algorithms scale as $\mathcal{O}(\chi^7)$ \cite{Lin2022}, offering a favorable alternative to the $\mathcal{O}(\chi^{10})$ scaling of many generic PEPS algorithms. Several isoPEPS algorithms have been recently developed for computing a single algebraically smallest eigenpair, i.e., the ground-state, and real-time evolution of nearest neighbor two-dimensional systems \cite{zaletel2020isometric, Haghshenas2019, Kadow2023, Wu2025, Sappler2025, dai2025fermionic}. However, methods for computing several of the algebraically smallest eigenpairs, i.e., excited states, with isoPEPS remain unexplored. 

\begin{figure*}[t]
\centering
\includegraphics[scale=1]{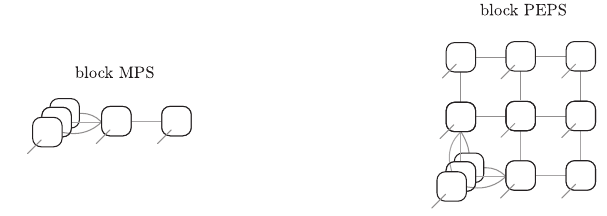}
\caption{Tensor network diagrams for representing three tensors. The block MPS network, shown here with 3 sites, provides a compact representation of states with limited entanglement in 1D geometries. In this paper we propose the block (isometric) PEPS ansatz, shown here with 9 sites, as a compact representation of several states with limited entanglement in 2D geometries. }
\label{fig:block_tns}
\end{figure*}

In this paper, we introduce a new block-isoPEPS tensor network ansatz shown in \cref{fig:block_tns} (right) for computing several eigenpairs of nearest neighbor two-dimensional quantum many-body systems. Our approach generalizes the state averaged MPS (or block tensor train) framework \cite{Pivzorn2012, dolgov2014, kressner2014, Li2024} \cref{fig:block_tns} (left), a well-established tensor network ansatz for computing several eigenpairs in one-dimensional systems, to two dimensions. The block-isoPEPS ansatz offers several desirable properties for handling nearest-neighbor Hamiltonians in two-dimensions. It retains the properties of isoPEPS, providing direct control over truncation errors through singular values and enabling efficient computation of local observables. Additionally, the isometric constraints enable exact and efficient orthogonalization of a block of tensor network states. To move the orthogonality center and block index around the tensor network we present the block Moses Move, an extension of the Moses Move introduced in \cite{zaletel2020isometric}. Using the block Moses Move we propose a new isoPEPS subspace iteration method for computing several eigenpairs with computational complexity scaling as $\bigO(\chi^7)$. Our approach differs from existing PEPS methods \cite{Vanderstraeten2015, vanderstraeten2019variational, Tu2024, zhang2025}, which do not impose orthogonality constraints during the optimization. One common strategy is to construct the excitation ansatz in the tangent space of an optimized ground state and then solve the effective Hamiltonian problem within this non-orthogonal subspace to extract excitation energies. Due to the restricted form of the ansatz, such approaches can only capture local excitations above the ground state. Alternatively, one can avoid these special excitation ansatz by directly optimizing the functional defined by the inverse norm matrix multiplied by the Hamiltonian matrix~\cite{zhang2025} within the variational space. 
We benchmark the proposed block-isoPEPS subspace iteration method on the two-dimensional transverse field Ising and Heisenberg models, demonstrating its accuracy, scalability, and ability to capture low-lying excitations across various parameter regimes. We also compare the proposed method with existing PEPS tangent space excitation methods and show that the proposed method is capable of capturing excitations that are not captured in the tangent space of the PEPS ground state with small bond dimension. 

The remainder of this paper is organized as follows. In \cref{sec:1disotns}, we review isometric tensor networks in one-dimension, i.e., the canonical forms of MPS, and block-MPS for representing multiple tensors. In \cref{sec:2d_isotns}, we introduce the block-isoPEPS ansatz which generalizes the block-MPS framework to two-dimensions. We also present the block Moses Move algorithm for moving the block site and orthogonality center through the network. In \cref{sec:subspace} we describe a subspace iteration algorithm for block-isoPEPS based on Trotter splitting of the matrix exponential. In \cref{sec:numerics} we demonstrate the block-isoPEPS subspace iteration on the transverse field Ising and Heisenberg models. We summarize our main findings and discuss future directions in \cref{sec:conclusions}.

\section{Isometric tensor networks in 1D} \label{sec:1disotns}


The operator \eqref{eq:Ham} acts on vectors in $\mathbb{C}^{d^L}$, which can be naturally identified with the tensor space $\mathbb{C}^{d \times \cdots \times d}$ through reshaping. We therefore represent states as tensors $T(\sigma_1,\ldots,\sigma_L)$ with indices $\sigma_j \in \{1,\ldots,d\}$. 
The number of entries in $T$ grows exponentially with the number of sites $L$ which makes computing and storing $T$ directly prohibitively expensive. Many tensors of interest can be computed and stored efficiently as a tensor network, which is a factorization of $T$ into low-dimensional tensors. 

\subsection{Matrix Product States (MPS)} \label{sec:block-mps}

A tensor $T$ is a matrix product state (MPS) if it is expressed as 
\begin{equation} \label{eq:MPS}
\begin{aligned}
    T(\sigma_1,\ldots,\sigma_L) &= A_1^{\sigma_1} A_2^{\sigma_2} \cdots A_L^{\sigma_L} \\
    &= \raisebox{-0.55\height}{\includegraphics[width=0.3\textwidth]{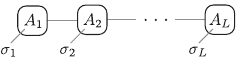}}
\end{aligned}
\end{equation}
Here, each $A_i^{\sigma_i}$ is a matrix associated with the index $\sigma_i$. For $i = 2,\ldots,L-1$, the $A_i^{\sigma_i}$ are $\chi \times \chi$ matrices, while the boundary tensors $A_1^{\sigma_1}$ and $A_L^{\sigma_L}$ are of size $1 \times \chi$ and $\chi \times 1$, respectively. The parameter $\chi$, known as the bond dimension, quantifies the amount of entanglement in the quantum state corresponding to the tensor $T$. 
In \eqref{eq:MPS} we have also used tensor diagram notation to the MPS where each tensor is depicted as a shape, and its indices (or dimensions) are shown as legs emanating from that shape. Connected legs between tensors represents a contraction over the shared indices. We have included labels for indices $\sigma_i$ in \eqref{eq:MPS}, however below we often do not label these indices as they can be inferred from the tensors they are attached to.

\subsection{MPS canonical forms}
The tensors in the MPS representation \eqref{eq:MPS} are not unique and many MPS algorithms rely on MPS tensors satisfying isometry constraints. The tensors in the MPS \eqref{eq:MPS} can be made isometries by performing QR or LQ decompositions on reshapings of the three-dimensional tensors $A_i$ into matrices, e.g., 
\begin{equation} \label{eq:mps_qr}
\begin{aligned}
    {\includegraphics[width=0.8\textwidth]{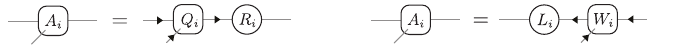}}
\end{aligned}
\end{equation}
In this diagrammatic notation isometries are indicated with arrows so that contracting a tensor with its conjugate over all indices with incoming arrows yields the identity operator in the space of all indices with outgoing arrows. Performing a sequence of QR decompositions and absorbing the $R_i$ into neighboring tensors on $A_1,\ldots,A_{j-1}$ and LQ decompositions and absorbing the $L_i$ into neighboring tensors on $A_{L}, \ldots,A_{j+1}$ the MPS \eqref{eq:MPS} is brought into the $j$th canonical (or isometric) form 
\begin{equation} \label{eq:mps_can}
\begin{aligned}
T(\sigma_1,\ldots,\sigma_L) 
= \raisebox{-0.5\height}
{\includegraphics[width=0.35\textwidth]{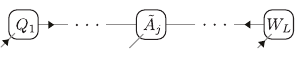}}
\end{aligned}
\end{equation}
where $A_j$ is referred to as the center of orthogonality. 
The computational cost of transforming a MPS into isometric form with QR-decomposition scales as $\mathcal{O}(dL\chi^3)$.

\subsection{Block MPS}
The MPS format in equation~\eqref{eq:MPS} representing a single state can be extended to represent a block of states by introducing an additional index on one tensor 
\begin{equation} \label{eq:block_MPS}
\begin{aligned}
T_{\alpha} &= A_1^{\sigma_1} \cdots A_{j,\alpha}^{\sigma_j} \cdots A_L^{\sigma_L} \\
&= \raisebox{-0.31\height}{\includegraphics[width=0.35\textwidth]{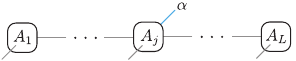}}
\end{aligned}
\end{equation}
where $\alpha = 1,2,\ldots,p$ indexes the block of $p$ states. It is convenient to work with block-MPS in canonical form where the center of orthogonality carries the block index $\alpha$. In this case, the inner product of two states $T_{\alpha}, T_{\beta}$ in the block is exactly the inner product of the low-dimensional tensors $A_{j,\alpha}^{\sigma_j}, A_{j,\beta}^{\sigma_j}$. This enables efficient and exact orthogonalization of the state block, e.g., by performing Gram-Schmidt on the $p$ MPS tensors $A_{j,\alpha}^{\sigma_j}$. 

The block index and orthogonality center can be moved through the MPS network using a similar QR-decomposition (or LQ-decomposition) in \eqref{eq:mps_qr} where the QR is performed so that the block index is passed to the $R$ factor 
\begin{equation} \label{eq:bmps_can}
\begin{aligned}
{\includegraphics[width=0.35\textwidth]{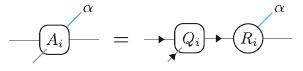}}
\end{aligned}
\end{equation}
Then the block index is passed to MPS tensor $i+1$ when $R_i$ is absorbed into its right neighbor $A_{i+1}$ 
\begin{equation} \label{eq:bmps_shift}
\begin{aligned}
{\includegraphics[width=0.35\textwidth]{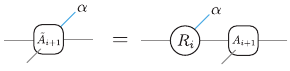}}
\end{aligned}
\end{equation}
The computational cost of this move scales as $\mathcal{O}(dLp\chi^3)$. 
In the algorithms presented in this paper, the block index $\alpha$ is never contracted over and hence it is not important that it has an isometry arrow. 


\section{Isometric tensor networks in 2D} \label{sec:2d_isotns}

Projected Entangled Pair State (PEPS) is a natural generalization of Matrix Product State (MPS) with additional connectivity between tensors enabling them to efficiently capture entanglement in quantum states defined on 2D geometries. To describe PEPS further, let us assume that the $d$-dimensional local vector spaces on which the $M_j$ in \eqref{eq:Ham} act are arranged in a rectangular lattice with $L_x$ sites in the vertical direction (number of rows) and $L_y$ sites in the horizontal direction (number of columns), for a total of $L = L_x \cdot L_y$ sites. Just like in \cref{sec:1disotns} we consider tensors except this time with the indices labeled to reflect the two-dimensional lattice geometry $T(\sigma_{11},\ldots,\sigma_{L_xL_y}) \in \mathbb{C}^{d \times \cdots \times d}$. 
A tensor $T$ is a PEPS if it is expressed as 
\begin{equation} \label{eq:peps}
\begin{aligned}
T(\sigma_{11},\ldots,\sigma_{L_xL_y}) = \mathcal{F}\left(A_{11}^{\sigma_{11}}, \ldots, A_{L_x L_y}^{\sigma_{L_x L_y}} \right) 
= \raisebox{-0.45\height}{\includegraphics[width=0.25\textwidth]{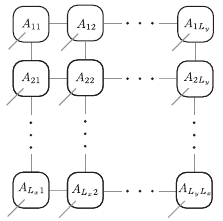}}
\end{aligned}
\end{equation}
where each interior tensor $A_{ij}$ has 5 legs and 4 neighbors. Due to the additional connectivity, it is not straightforward to write the PEPS as a product of matrices like with MPS in \eqref{eq:MPS}. Instead, we let $\mathcal{F}$ denote the tensor contraction over the connected legs in the tensor network diagram shown in \eqref{eq:peps}. As mentioned in \cref{sec:intro}, a computational challenge of PEPS is that the cost of computing inner products exactly, such as those required for observables, scales exponentially with the system size. To circumvent this challenge, we restrict to a subclass of PEPS where most tensors in the network are isometries as was done in \cite{zaletel2020isometric}. This approach can be considered a generalization of MPS canonical form \eqref{eq:mps_can} to PEPS. However unlike MPS, the orthogonality center of PEPS cannot be moved around the network without introducing errors or increasing bond-dimension. To move the orthogonality center in PEPS one can decompose a column into the product of an isometric column and a non-isometric column without physical indices using the so-called Moses Move 
\begin{equation} \label{eq:peps_qr}
\begin{aligned}    {\includegraphics[width=0.3\textwidth]{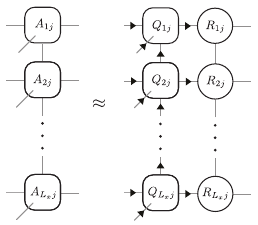}}
\end{aligned}
\end{equation}
which plays the role of the QR-decomposition \eqref{eq:mps_qr} used to move the orthogonality center in MPS. Variational and greedy algorithms for performing the approximate decomposition were proposed in \cite{zaletel2020isometric}. After obtaining the column decomposition \eqref{eq:peps_qr} the non-isometric column can be absorbed into the neighbor column. By grouping indices appropriately this product can be viewed as a product of matrix product operators (MPOs) for which several algorithms are available. The Moses Move can also be applied to decompose rows of a PEPS. After performing a suitable sequence of Moses Moves, any PEPS can be approximated by an isometric PEPS with orthogonality center in the bottom left corner
\begin{equation} \label{eq:peps_can}
\begin{aligned}
T(\sigma_{11},\ldots,\sigma_{LL}) 
\approx \raisebox{-0.5\height}
{\includegraphics[width=0.25\textwidth]{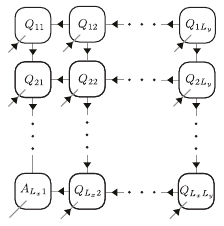}}
\end{aligned}
\end{equation} 
We refer to PEPS with a center of orthogonality such as \eqref{eq:peps_can} as isometric PEPS (isoPEPS). It is important to note that unlike in MPS, the orthogonality center of a isoPEPS cannot be shifted exactly without increasing the bond dimension. This implies that the set of isoPEPS with a fixed bond dimension forms a proper subset of generic PEPS. The error introduced when moving the orthogonality center depends on the specific algorithm used, such as the sequential Moses Move \cite{zaletel2020isometric}. While variational algorithms have also been proposed for this task, the greedy approach of the sequential Moses Move has been shown to achieve comparable accuracy with a lower computational cost. 



\subsection{Block isometric PEPS (block-isoPEPS)}

Next we introduce block-isoPEPS, an extension of isoPEPS for representing several tensors. Similar to block-MPS in \eqref{eq:block_MPS}, block-isoPEPS adds an additional block index $\alpha$ to the tensor at the center of orthogonality to represent several tensors that share the surrounding isometric tensors 
\begin{equation} \label{eq:bpeps}
\begin{aligned}
T_{\alpha}(\sigma_{11},\ldots,\sigma_{L_xL_y}) &= \mathcal{F}\left(Q_{11}^{\sigma_{11}}, \ldots, A_{ij,\alpha}^{\sigma_{ij}}, \ldots, Q_{L_x L_y}^{\sigma_{L_x L_y}} \right) 
= \raisebox{-0.45\height}{\includegraphics[width=0.25\textwidth]{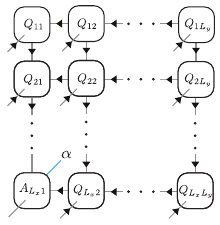}}
\end{aligned}
\end{equation}
where $A_{ij,\alpha}$ (in this case $i=L_x$ and $j=1$) carries the block index $\alpha = 1, \ldots, p$. To enable optimization and time evolution with block-isoPEPS, we extend the algorithms for shifting the orthogonality center introduced in \cite{zaletel2020isometric} to also move the block index around the isoPEPS network. We adopt the greedy strategy known as the sequential Moses Move \cite{zaletel2020isometric, Lin2022, dai2025fermionic}. 

\subsubsection{Block sequential Moses Move} \label{sec:block_mm}

We describe the block sequential Moses Move on the $j$th column $C_j$ of an isoPEPS assuming the orthogonality center and block index are positioned at the bottom of the column 
\begin{equation} \label{eq:col}
    C_j = \raisebox{-0.45\height}{\includegraphics[width=0.1\textwidth]{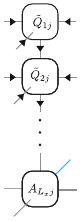}}
\end{equation}
The move begins by decomposing the bottom tensor $A_{L_xj}$ into a tensor ring \cite{TensorRing2020} consisting of three tensors with isometric constraints 
\begin{equation} \label{eq:tr1}
{\includegraphics[width=0.3\textwidth]{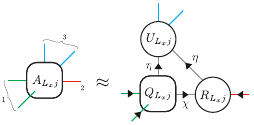}}
\end{equation}
where indices are grouped into three dimensions as shown. The user prescribes a maximum vertical bond-dimension $\eta$ and a maximum horizontal bond-dimension $\chi$ for the tensor ring decomposition \eqref{eq:tr1}. Details on computing the isometric tensor ring decomposition are provided in \cref{sec:iso_tr}. Then the $U_{L_x j}$ tensor is absorbed into its above neighbor $\tilde{Q}_{(L_x-1) j}$ which yields 
\begin{equation*} \label{eq:Up1}
{\includegraphics[width=0.3\textwidth]{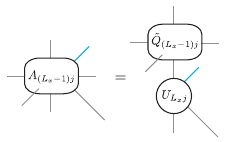}}
\end{equation*}
From here we proceed recursively for $i=L_x-1, \ldots, 2$ performing an isometric tensor ring decomposition of $A_{ij}$ 
\begin{equation} \label{eq:tr2}
{\includegraphics[width=0.3\textwidth]{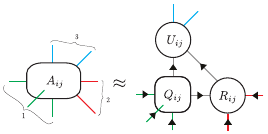}}
\end{equation}
once again with maximum vertical bond-dimension $\eta$ and horizontal bond-dimension $\chi$, and then absorbing the $U_{ij}$ tensor to its upper neighbor $\tilde{Q}_{i-1j}$ to form 
\begin{equation*} \label{eq:Up2}
{\includegraphics[width=0.3\textwidth]{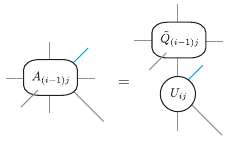}}
\end{equation*}
Once we have reached $i=2$ and formed $A_{11}$, we decompose $A_{11}$ using a truncated SVD 
\begin{equation} \label{eq:last_svd}
{\includegraphics[width=0.3\textwidth]{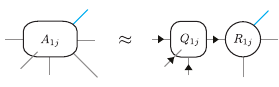}}
\end{equation}
with maximum horizontal bond-dimension $\chi$, where $R_{1j}$ absorbs the singular values. The end result of the block sequential Moses Move is an approximation of the original column \eqref{eq:col} as a product of an isometric column and a non-isometric column without physical indices just as with the sequential Moses Move \eqref{eq:peps_qr}
\begin{equation} \label{eq:bpeps_qr}
C_j = \raisebox{-0.51\height}{\includegraphics[width=0.3\textwidth]{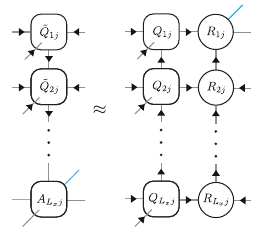}}
\end{equation}
where the vertical and horizontal bond-dimensions are less than $\eta$ and $\chi$, respectively. In addition, the block index has been shifted to the $R_{11}$ tensor. The computational cost of the tensor ring decomposition \eqref{eq:tr2} scales as $\bigO(\chi^3\eta^4d^2p + \chi^2\eta^5p^2)$ (see \cref{sec:iso_tr} for cost of isometric tensor ring), which dominates one step of the block sequential Moses Move. We perform $L_x$ of these steps for a total cost of $\bigO(L_x(\chi^3\eta^4d^2p + \chi^2\eta^5p^2))$ to split column $C_j$ in \eqref{eq:bpeps_qr}. 

Once the column $C_j$ has been split as in \eqref{eq:bpeps_qr}, the column containing the $R_{ij}$ tensors are contracted into the neighboring column $C_{j+1}$ making the top core of the $j+1$ column the orthogonality center with the block index 
\begin{equation} \label{eq:zipup}
C_{j+1} = \raisebox{-0.51\height}{\includegraphics[width=0.3\textwidth]{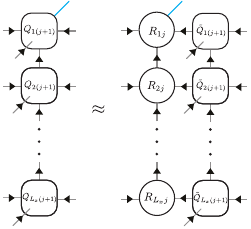}}
\end{equation}
During this tensor contraction the horizontal bond-dimension is truncated so that it is no larger than $\eta$. There are several possible algorithms one can employ to perform this contraction and truncation. 
If the contraction of $R_{ij}$ into its neighboring column is done naively, then the vertical bond dimension of the new column will be equal to the product of the vertical bond dimension of $R_{ij}$ and the vertical bond dimension of the neighboring column before contracting. The SVD truncation following such naive contraction costs  
$\bigO(\chi^5\eta^3p)$.
This cost can be reduced to 
$\bigO(\chi^4\eta^3p)$ 
by adapting techniques designed for MPO-MPS multiplication such as variational method, zip-up method \cite{Stoudenmire2010}, or randomized SVD, which limit bond-dimension growth by performing approximations during the contraction procedure. 

In our numerical experiments we employ the zip-up method which controls the growth of the vertical bond-dimension during contraction. Because the tensor network is not maintained in canonical form during this process, these local truncations are not globally optimal. To mitigate the resulting uncontrolled truncation errors, we first orthogonalize the $R$ column from top to bottom, thereby placing the isometric center at the bottom of the $R$ column. In our experience the zip-up algorithm with an appropriately chosen singular value truncation threshold provides a substantial computational speedup compared to direct contraction methods without introducing significant error. The total cost of shifting the orthogonality center and block tensor with block-sequential Moses Move and zip-up method is $\bigO(L_x(\chi^3\eta^4d^2p + \chi^2\eta^5p^2))$. 

\subsubsection{Isometric tensor ring decomposition} \label{sec:iso_tr}

A key subroutine and main source of error in the (block) sequential Moses Move is the decomposition of a three-dimensional tensor 
\begin{equation*}
    {\includegraphics[width=0.1\textwidth]{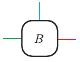}}
\end{equation*}
into a tensor ring \cite{TensorRing2020} with isometric constraints and maximum vertical and horizontal bond-dimensions $\eta$ and $\chi$ as in \eqref{eq:tr1} and \eqref{eq:tr2}. We follow the greedy method proposed in \cite{zaletel2020isometric}, which uses a sequence of two truncated SVDs. In between these two SVDs, one can optimize a unitary disentangler \cite{Wei2025} to reduce the error in the second truncated SVD. 

The first SVD is used to separate the left index from the up and right indices 
\begin{equation} \label{eq:svd_1}
{\includegraphics[width=0.3\textwidth]{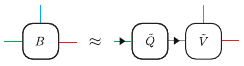}}
\end{equation}
with the singular values absorbed into $\tilde{V}$. The maximum bond-dimension of the new bond is $\eta \chi$. Then the new bond is split into two bonds, one with maximum dimension $\eta$ and one with maximum dimension $\chi$ by reshaping $\tilde{Q}$ and $\tilde{V}$ 
\begin{equation} \label{eq:split_bond}
{\includegraphics[width=0.4\textwidth]{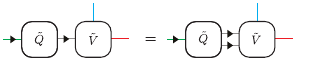}}
\end{equation}
At this point the isometric tensor ring can be completed by performing a truncated SVD on $\tilde{V}$ that separates the bonds on either side of the dashed line below 
\begin{equation} \label{eq:svd2}
{\includegraphics[width=0.3\textwidth]{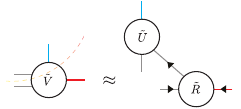}}
\end{equation}
using maximum bond-dimension $\eta$ 
with the singular values absorbed into $\tilde{U}$. 
However instead of performing an SVD on $\tilde{V}$ the truncation error can be reduced by optimizing a unitary operator $D$ referred to as a disentangler to obtain 
\begin{equation} \label{eq:dis}
{\includegraphics[width=0.35\textwidth]{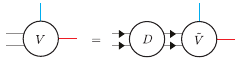}}
\end{equation}
where the objective of optimizing $D$ is to reduce the error of the truncated SVD 
\begin{equation} \label{eq:svd3}
{\includegraphics[width=0.3\textwidth]{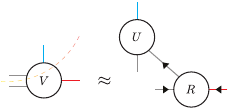}}
\end{equation}
For more details on the numerical optimization of disentanglers see \cite{Wei2025}. Finally setting 
\begin{equation}
{\includegraphics[width=0.3\textwidth]{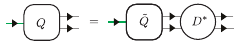}}
\end{equation}
completes the isometric tensor ring decomposition 
\begin{equation}
{\includegraphics[width=0.3\textwidth]{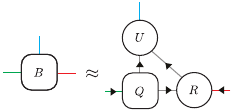}}
\end{equation}

To improve the accuracy of the isometric tensor ring decomposition one can proceed to optimize the tensor cores several times using an alternating least squares (ALS) approach. However in our experience we found that ALS does not significantly improve the accuracy of the greedy approach using two SVDs and a disentangler. 

To estimate the computational cost of the greedy isometric tensor ring decomposition algorithm we assume the dimension of the left, up, right indices of $B$ are $n_1, n_2, n_3$, respectively. The first SVD \eqref{eq:svd_1} costs 
$\bigO(n_1^2n_2n_3)$, assuming $n_1^2 \leq (n_2n_3)^2$. Optimizing the disentangler \eqref{eq:dis} and computing the final SVD \eqref{eq:svd3} costs $\bigO((\eta\chi)^3 + (\eta \chi)^2 n_2 n_3 + (\eta n_2)^2 \chi n_3)$. 




\section{Subspace iteration for nearest-neighbor Hamiltonians in 2D} \label{sec:subspace}

We now present an inexact subspace iteration based on block-isoPEPS for \eqref{eq:Ham} with nearest-neighbor interactions on an $L_x \times L_y$ lattice. Using two-dimensional site indices, \eqref{eq:Ham} can be rewritten as 
\begin{equation} \label{eq:Ham2}
H = \sum_{i=1}^{L_x-1} \sum_{j=1}^{L_y} H_{(i,j),(i+1,j)} + 
\sum_{i=1}^{L_x} \sum_{j=1}^{L_y-1} H_{(i,j),(i,j+1)} 
\end{equation}
where $H_{(i,j),(i+1,j)}$ acts on the vertical pair of sites $(i,j)$ and $(i+1,j)$, and $H_{(i,j),(i,j+1)}$ acts on the horizontal pair $(i,j)$ and $(i,j+1)$. 


A standard approach for computing the eigenvector corresponding to the algebraically smallest eigenvalue of \eqref{eq:Ham} is to apply the power method to the matrix exponential $e^{-H}$. Directly evaluating the full matrix exponential is computationally intractable. Instead, we consider $e^{-\tau H}$ for a small time step $\tau>0$, which is more amenable to computation due to the structure of \eqref{eq:Ham2}. Using the Suzuki–Trotter decomposition we obtain a first order approximation 
\begin{equation} \label{eq:lie_trotter}
\begin{aligned}
e^{-\tau H} &\approx \prod_{i=1}^{L_x-1} \prod_{j=1}^{L_y} e^{-\tau H_{(i,j),(i+1,j)}}
\prod_{i=1}^{L_x} \prod_{j=1}^{L_y-1} e^{-\tau H_{(i,j),(i,j+1)}}, 
\end{aligned}
\end{equation} 
that can be conveniently applied to a block-isoPEPS. Repeated application of \eqref{eq:lie_trotter} to a starting vector is an inexact power iteration referred to as imaginary time evolution in the physics literature. We apply \eqref{eq:lie_trotter} to evolve a subspace spanned by $p$ basis tensors represented as a block-isoPEPS \eqref{eq:bpeps}. To ensure the basis tensors converge to different eigenvectors we orthogonalize the tensors in the block after each application of the split matrix exponential \eqref{eq:lie_trotter}. 

\subsection{Applying the split matrix exponential} 
\label{sec:matvec}
We follow the \TEBD2 algorithm proposed in \cite{zaletel2020isometric}.
Each factor in the product \eqref{eq:lie_trotter} is a $d^L \times d^L$ matrix that acts non-trivially on only two sites. Its action on a tensor can be computed by contracting a $d \times d \times d \times d$ gate tensor to the corresponding pair of tensor indices. To simplify our presentation we assume that all gate tensors are identical and denote them by $G$, although in practice they can differ. To apply the first-order approximation \eqref{eq:lie_trotter} to a block-isoPEPS $T_{\alpha}$ we assume that the orthogonality center and block index are located at site $(1,1)$ 
\begin{equation} \label{eq:bpeps2}
\begin{aligned}
T_{\alpha} = \raisebox{-0.45\height}{\includegraphics[width=0.25\textwidth]{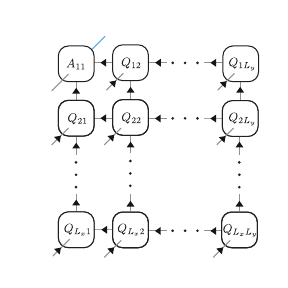}}
\end{aligned}
\end{equation}
We apply the vertical nearest-neighbor gates recursively for $j=1,2,\ldots,L_y$. To apply the gates in column $j$ we proceed recursively over the rows $i=1,2,\ldots,L_x-1$. For each of these nearest-neighbor gates we perform the tensor contraction 
\begin{equation} \label{eq:gate}
{\includegraphics[width=0.2\textwidth]{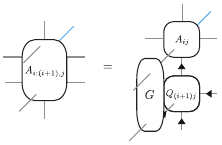}}
\end{equation}
Then we reshape the resulting tensor $A_{i:(i+1),j}$ so that the block index is grouped with the bottom indices and perform a truncated SVD with maximum bond-dimension $\eta$ to separate the indices on either side of the dashed line 
\begin{equation} \label{eq:svd_tebd}
{\includegraphics[width=0.2\textwidth]{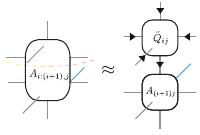}}
\end{equation}
absorbing the singular values into $A_{(i+1)1}$. Repeating this procedure for $i = 1, \ldots, L_x - 1$ produces an updated column $C_j$ in the isometric form shown in \eqref{eq:col} with all relevant vertical nearest-neighbor gates applied. Next we use the block-sequential Moses Move \eqref{eq:col}-\eqref{eq:zipup} to shift the orthogonality center and block index to the top of the $j+1$, preparing $C_{j+1}$ for the application of the column $j+1$ vertical gates. This process is repeated recursively for $j=1,\ldots,L_y$ at which point all vertical gates have been applied. After applying the vertical gates to the final column, we move the orthogonality center to the top of the isoPEPS using a sequence of QR-decompositions. To apply the horizontal gates we rotate the block-isoPEPS counter-clockwise by ninety degrees and then reuse the above implementation on columns.  The computational cost of applying a vertical or horizontal gate followed by truncation is dominated by the SVD in \eqref{eq:svd_tebd} which has computational cost scaling as $\bigO(\chi^6 \eta^3 d^3 p)$ for interior PEPS tensors. 

This cost can be significantly reduced using a technique known as the reduced update \cite{Corboz2010, lubasch2014algorithms}. The idea is to first decrease the size of the tensors to which the gate is applied by performing QR decompositions 
\begin{equation} \label{eq:reduced_cores}
{\includegraphics[width=0.2\textwidth]{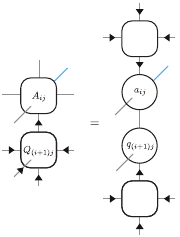}}
\end{equation}
and then applying the two-site gate and truncated SVD to the smaller tensors 
\begin{equation} \label{eq:reduced_svd}
{\includegraphics[width=0.2\textwidth]{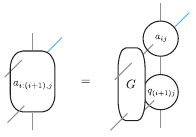}}, 
\quad 
\raisebox{-0.14\height}
{\includegraphics[width=0.3\textwidth]{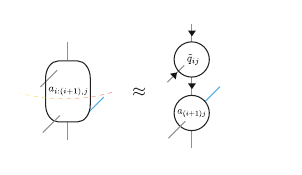}}
\end{equation} 
Assuming the horizontal bond-dimensions are $\chi$ and vertical bond-dimensions are $\eta$, the computational cost of QR-decomposition for the upper core carrying the block index is $\bigO(\chi^4\eta^3 d p)$ and for the bottom core the cost is $\bigO(\chi^4\eta^3 d)$. These QR-decompositions dominate the cost of the reduced update. 

\subsection{Orthogonalizing the tensor block}

To ensure that the tensors represented by the block-isoPEPS \eqref{eq:bpeps} converge to distinct eigenvectors, we orthogonalize the tensors in the block after each application of the split matrix exponential. Owing to the isometric structure of the PEPS tensors, the inner product between two tensors $T_{\alpha}$ and $T_{\beta}$ can be computed efficiently and exactly from the inner products of the corresponding low-dimensional tensors $A_{ij,\alpha}$ and $A_{ij,\beta}$. Consequently, orthogonalizing the block amounts to orthogonalizing the small tensors $A_{ij,\alpha}$ at the center of orthogonality in the block-isoPEPS. This is accomplished using the Gram–Schmidt procedure.

\subsection{Calculating eigenvalues and observables} 

The subspace iteration described above yields approximate eigenvectors of \eqref{eq:Ham2}, however it does not produce approximations to the eigenvalues. To compute approximate eigenvalues from a block-isoPEPS approximation $T_{\alpha}$ of a block of eigenvectors, we compute approximations of the Rayleigh quotients 
\begin{equation} \label{eq:RQ}
\begin{aligned}
E_{\alpha} 
&= 
\frac{\left\langle T_{\alpha}, H T_{\alpha}\right\rangle}
{\|T_{\alpha}\|} \\
&= 
\sum_{i=1}^{L_x-1} \sum_{j=1}^{L_y} 
\frac{\left\langle T_{\alpha}, H_{(i,j),(i+1,j)} T_{\alpha}\right\rangle}
{\|T_{\alpha}\|}
+ 
\sum_{i=1}^{L_x} \sum_{j=1}^{L_y-1} 
\frac{\left\langle T_{\alpha}, H_{(i,j),(i,j+1)} T_{\alpha}\right\rangle}
{\|T_{\alpha}\|}. 
\end{aligned}
\end{equation}
To accomplish this, we follow exactly the same procedure described in \cref{sec:matvec}, however instead of applying the gate tensor $G$ as in \eqref{eq:gate}, we compute one of the terms in the sum \eqref{eq:RQ} by performing the tensor contraction
\begin{equation}
\begin{aligned}
\left\langle T_{\alpha}, H_{(i,1),(i+1,1)} T_{\alpha}\right\rangle =
\raisebox{-0.45\height} {\includegraphics[width=0.2\textwidth]{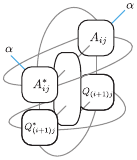}}, \qquad
\|T_{\alpha}\|^2 = \raisebox{-0.45\height} {\includegraphics[width=0.2\textwidth]{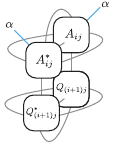}}
\end{aligned}
\end{equation}
Other observable can be computed in a similar fashion. 

\subsection{Errors and computational cost}
\label{sec:errors_and_cost}

In all block-isoPEPS algorithms discussed above, truncated SVDs are applied to tensors at the orthogonality center of the PEPS. This design ensures that the error introduced by each truncation is directly controlled by the discarded singular values. Note that this is not the case when truncating bonds in a general PEPS. In the subspace iteration algorithm the sources of truncation error are the SVD after applying the two-site gates \eqref{eq:gate}, and the two truncated SVDs in the block-sequential Moses Move \eqref{eq:svd_1} and \eqref{eq:svd3} used to construct the isometric tensor ring decomposition. There is also $\bigO(\tau^2)$ error due to the splitting of the matrix exponential \eqref{eq:lie_trotter}. Note that this splitting error can be reduced by using higher order splitting algorithms for the matrix exponential. 

The cost of a single iteration of block-isoPEPS subspace iteration is dominated by applying the split matrix exponential. To estimate the computational cost we assume a square lattice $L_x=L_y$. 
Application of each $(L_x-1)^2$ gate followed by truncation costs $\bigO(\chi^4\eta^3 d p)$ when using the reduced update \eqref{eq:reduced_cores}-\eqref{eq:reduced_svd} for a total cost of $\bigO(L_x^2 \chi^4\eta^3 d p)$. 
In addition we perform $(L_x-1)^2$ block sequential Moses Moves each of which cost $\bigO(L_x(\chi^3\eta^4d^2p + \chi^2\eta^5p^2))$ for a total cost of $\bigO(L_x^2(\chi^3\eta^4d^2p + \chi^2\eta^5p^2))$. 
The computational cost of applying the split matrix exponential is thus $\bigO(L_x^2(\chi^4\eta^3 d p + \chi^3\eta^4d^2p + \chi^2\eta^5p^2))$, which scales as the seventh power of the bond dimension, $\bigO(\chi^4\eta^3 + \chi^3\eta^4 + \chi^2 \eta^5)$. This scaling with respect to bond-dimension is substantially lower than the $\bigO(\chi^{10})$ scaling of excited-state algorithms based on generic PEPS, such as the tangent-space method summarized in \ref{sec:peps_tangent_space}. However because isoPEPS is less expressive, it may require a larger bond dimension than generic PEPS to represent the same tensor. 



\section{Numerical experiments} \label{sec:numerics}

In this section, we present numerical demonstrations of the proposed block-isoPEPS subspace iteration algorithm. 
We computed ground and low-lying excited states for the transverse-field Ising (TFI) model
\begin{equation} \label{eq:TFI}
H_{\text{TFI}} = - \sum_{\langle i, j \rangle} \sigma_i^z \sigma_j^z - g \sum_i \sigma_i^x,
\end{equation}
and the Heisenberg model
\begin{equation} \label{eq:heis}
H_{\text{heis}} = \sum_{\langle i, j \rangle} \sigma_i^x \sigma_j^x + \sigma_i^y \sigma_j^y + \sigma_i^z \sigma_j^z,
\end{equation}
on a two-dimensional square lattice of side length $L_x=L_y$ ($L=L_x^2$ total sites) with open boundary conditions. The operators $\sigma_i^a$ ($a \in {x, y, z}$) denote Pauli matrices acting on site $i$, and $\langle i, j \rangle$ indicates nearest-neighbor pairs. 

To benchmark the block-isoPEPS method, reference eigenvalues $E_{\alpha,\text{ref}}$ for the TFI model \eqref{eq:TFI} were computed using the block MPS alternating optimization algorithm described in \cite{dolgov2014}. The index $\alpha=0$ corresponds to the ground state (the algebraically smallest eigenvalue), while $\alpha=1,2,\ldots$ label the first several excited states. With a maximum bond dimension of $5000$, reference energies for the TFI model were obtained with relative accuracy of order $10^{-6}$ for $L_x \leq 8$ and $10^{-4}$ for $L_x > 8$. 
We also considered the antiferromagnetic Heisenberg model, with reference energies obtained using quantum number–conserving DMRG in ITensor \cite{ITensor}. The ground state belongs to the $S_z=0$ symmetry sector and the first excited state belongs to the $S_z=1$ sector. We used a maximum MPS bond-dimension of $1000$ and the DMRG truncation error was on the order of $10^{-5}$. 

\begin{figure}[t]
\centering
\includegraphics[scale=.5]{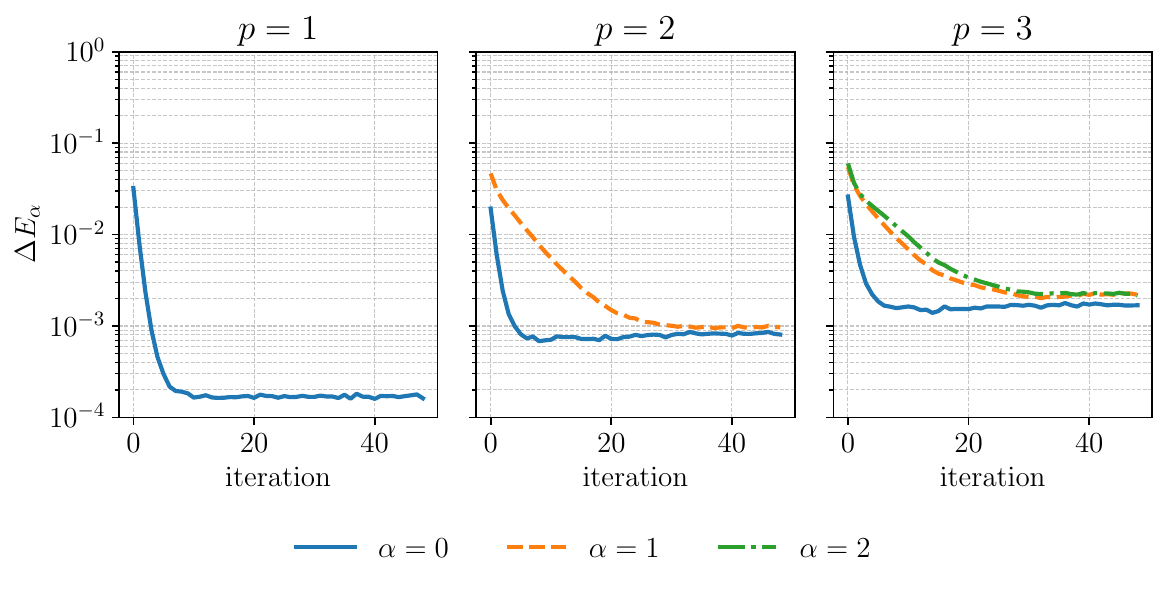}
\caption{
Relative eigenvalue error for the 2D transverse field Ising (TFI) model with $g=3.5$ on a $6 \times 6$ lattice. Eigenvalues were computed using the block-isoPEPS subspace iteration algorithm with bond dimensions $\chi = 8$ and $\eta = 16$. Results are shown for block sizes $p = 1$, $2$, and $3$. } 
\label{fig:L6_Ising_block_sizes}
\end{figure}

\begin{figure*}[t]
\centering

\begin{minipage}{0.48\textwidth}
\centering
{\footnotesize $g = 3.5$}\\
\includegraphics[scale=0.4]{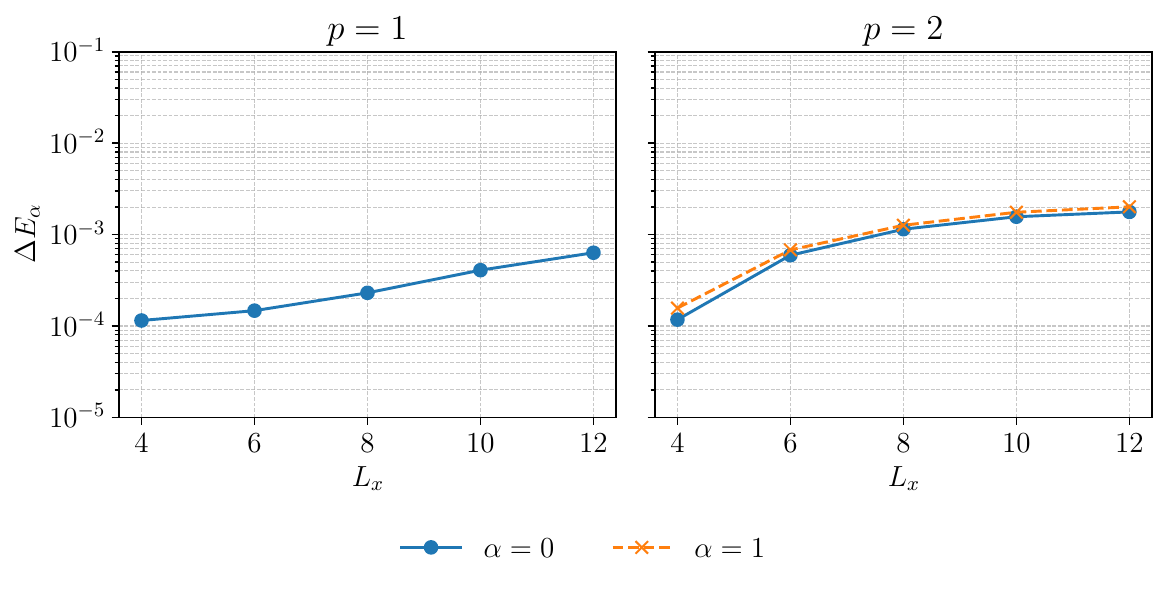}
\end{minipage}
\hfill
\begin{minipage}{0.48\textwidth}
\centering
{\footnotesize $g = 3.0$}\\
\includegraphics[scale=0.4]{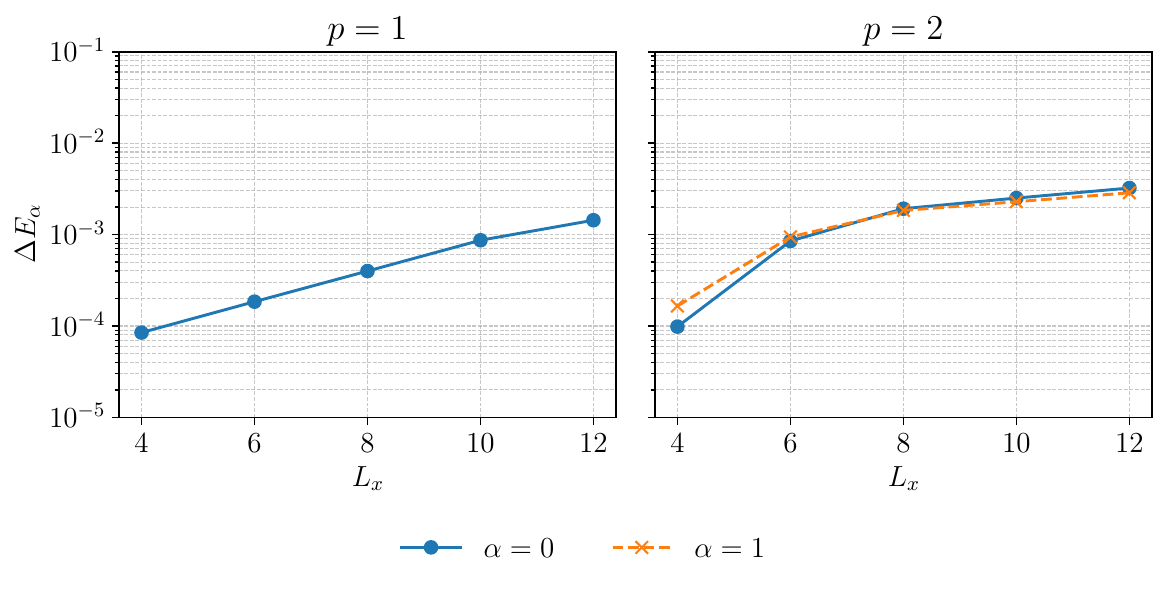}
\end{minipage}

\caption{
Relative eigenvalue error for the transverse field Ising (TFI) model \eqref{eq:TFI} with $g=3.5$ and $g=3.0$ as a function of lattice side length $L_x$. Approximate eigenvalues were computed using $50$ iterations of block iso-PEPS subspace iteration with bond dimensions $\chi = 12$ and $\eta = 20$ and $\tau = 0.1$. Results are shown for block sizes $p = 1$ and $p = 2$.}
\label{fig:Ising_p2_versus_L}
\end{figure*}

To move the orthogonality center within the block-isoPEPS, we employed the block-sequential Moses Move described in \cref{sec:block_mm}. For the column containing the orthogonality center, we used a maximum bond dimension $\eta$, while the surrounding isometric tensors used a bond dimension $\chi \leq \eta$, following \cite{zaletel2020isometric, Kadow2023}. To compute a disentangler \cite{Wei2025} for each isometric tensor ring decomposition we minimized the Reny-1/2 entropy using a Riemannian Conjugate Gradient method with a maximum of $30$ iterations. After performing the block-sequential Moses Move, we contracted the $R$ column into its neighbor using the zip-up algorithm with a singular value tolerance of $10^{-6}$. To evaluate the accuracy of the block-isoPEPS approximation, we used the relative eigenvalue error 
\begin{equation} \label{eq:rel_eig_err}
\Delta E_{\alpha} = \left|\frac{E_{\alpha} - E_{\alpha,\text{ref}}}{E_{\alpha,\text{ref}}}\right|, 
\end{equation}
where $E_{\alpha}$ is the approximate energy computed with the proposed block-isoPEPS method. 

We first studied the 2D TFI model \eqref{eq:TFI} with $g=3.5$ on a $6 \times 6$ lattice. In this case, the first three eigenvalues corresponding to $\alpha=0,1,2$, i.e., the ground state and first two excited state energies, are non-degenerate. We ran the isoPEPS subspace iteration with $\tau = 0.1$, bond dimensions $\chi = 8$ and $\eta = 16$, and varied the block size $p = 1, 2, 3$. We compared the computed eigenvalues with reference results from block-MPS DMRG and display the relative error \eqref{eq:rel_eig_err} of the eigenvalues in \cref{fig:L6_Ising_block_sizes} versus iteration. We observe that the approximation error increases with block size $p$, which is expected as more vectors are represented using a block isometric tensor network. Increasing the block size from $p = 1$ to $p = 2$ increases the final error by nearly an order of magnitude, from approximately $10^{-4}$ to $10^{-3}$. Further increasing the block size $p = 2$ to $p = 3$ results in a smaller increase, with the errors for $p = 3$ remaining close to $10^{-3}$. We also observe in \cref{fig:L6_Ising_block_sizes} that larger block sizes can require more iterations to reach convergence. 

Next, we fixed the bond dimensions to $\chi = 12$, $\eta = 20$, and considered block sizes $p = 1, 2$. We then applied the isoPEPS subspace iteration algorithm with time step $\tau = 0.1$ to compute the ground state and first excited state of \eqref{eq:TFI} with $g = 3.5$ on lattices of increasing size $L_x = 4, 6, 8, 10, 12$. 
We ran 50 iterations for each simulation and in \cref{tab:cpu_timings} report the average CPU-time for a iteration. 
\begin{table}[t]
\centering
\begin{tabular}{|l|cc|}
\hline
{ } & $p=1$ & $p=2$ \\
\hline
$L_x=4$  & 00:02 & 00:05 \\
$L_x=6$  & 00:38 & 01:44 \\
$L_x=8$  & 03:29 & 09:32 \\
$L_x=10$ & 10:28 & 25:22 \\
$L_x=12$ & 21:49 & 46:40 \\
\hline
\end{tabular}
\caption{Average CPU time (mm:ss) for a single iteration of isoPEPS subspace iteration applied to the TFI model with $g=3.5$, using bond dimensions $(\chi,\eta) = (12,20)$ and block sizes $p = 1, 2$ on square lattices with side lengths $L_x = 4, 6, 8, 10, 12$.} 
\label{tab:cpu_timings}
\end{table}
We ran the same experiment with $g=3.0$, closer to the critical point $g_{c} \approx 3.04$ of the TFI \eqref{eq:TFI}. In \cref{fig:Ising_p2_versus_L}, we display the relative errors \eqref{eq:rel_eig_err} of the first two eigenvalues as a function of $L_x$ for both $g=3.5$ and $g=3.0$. The left panel shows the error scaling for $p = 1$, corresponding to an isoPEPS representing only the ground state. The right panel shows the error for $p = 2$, corresponding to a block-isoPEPS representing both the ground and first excited state. We observe that from $L_x = 4$ to $L_x = 6$, the error of the $p=2$ block-isoPEPS increases more than that of the $p = 1$ isoPEPS. However, from $L_x = 6$ to $L_x = 8$, the errors for both $p = 1$ and $p = 2$ increase at a similar rate and from $L_x=8$ to $L_x=12$ the error of $p=1$ isoPEPS increases at a faster rate than $p=2$ isoPEPS. Overall we find that, for both $g=3.5$ and $g=3.0$, computing the first excited state of the TFI model does not dramatically impact the error when directly compared to the isoPEPS computation of the ground state only. 

To further assess the isoPEPS results in relation to existing PEPS-based approaches, we applied the tangent-space method to finite PEPS for the $L_x=4$ TFI model across different phases. More details on the PEPS tangent-space method are provided in \ref{sec:peps_tangent_space}. The comparison between isoPEPS with block index $2$ and the PEPS tangent-space method is summarized in Table~\ref{tab:peps_isopeps_L4}. We consider different values of the transverse field $g$ to cover both the polarized phase ($g > g_c$) and the symmetry-broken phase ($g < g_c$). For isoPEPS with $\chi=4$, the errors remain uniform at the level of $10^{-3}$ across both phases, with higher accuracy deep in the symmetry-broken regime (e.g., $g=1.0$). The accuracy can be systematically improved by at least one order of magnitude when the bond dimension is increased to $\chi=12$. 

In contrast, the PEPS tangent-space method with $\chi=2$ achieves better accuracy in the polarized phase (e.g., $g=3.5$) but performs worse in the symmetry-broken regime. Starting from the ground state, we find that the ansatz fails to capture the first excited state at $\chi=2$ for $g=1.0$ and it performs worse in describing the first excitation at $g = 2.0$. This can be understood from the structure of the ansatz: in the symmetry-broken phase, the first excited state becomes degenerate with the ground state in the thermodynamic limit. In a finite system, this manifests as a small gap that decreases with $g$, as indicated by the reference energies listed in Table~\ref{tab:peps_isopeps_L4}. For such small energy gaps, a PEPS ansatz with limited bond dimension can approximate the ground state only as a superposition of the true ground state and the first excited state. Consequently, local excitations constructed on top of this approximate, symmetry-broken ground state cannot reach the other symmetry sector and therefore fail to reproduce the first excited state. This limitation can be mitigated by increasing $\chi$, with $\chi = 4$, the description of the ground state and the first excited state becomes equally accurate at $g = 2$. In contrast, at $g = 1$, one of the states is still not recovered.

This behavior stands in contrast to the block-isoPEPS method. In block-isoPEPS, the environment of the local tensor with block index is optimized simultaneously for both the ground state and the excited states, whereas in the tangent-space method the environment is fixed to that of the ground-state tensor. Consequently, nearly degenerate states are captured more accurately in block-isoPEPS, highlighting a potential advantage of employing block environments. For those states that can be described by the tangent-space ansatz, we further observe that increasing $\chi$ yields larger improvements in accuracy for PEPS than increasing $\chi,\eta$ in isoPEPS. However, the computational cost of PEPS grows more rapidly, scaling as $\bigO(\chi^{14})$ or $\bigO(\chi^{10})$ when restricting to a subspace of the tangent space, compared to isoPEPS, which scales as $\bigO(\chi^{7})$. A more detailed study of the trade-off between accuracy and computational cost for PEPS and isometric PEPS methods is left for future investigation. 


\begin{table*}[t]
\centering
\begin{tabular}{|c|c|ccccc|}
\hline
\textbf{$g$} & $\alpha$ & reference & isoPEPS $\chi=12,\eta=20$ & isoPEPS $\chi=4,\eta=8$ & PEPS $\chi=2$ & PEPS $\chi=4$ \\
\hline
\multirow{2}{*}{3.5} 
  & 0 & $-57.82436977$ & $1.17 \times 10^{-4}$ & $9.85 \times 10^{-4}$ & $9.68\times 10^{-6}$ & $2.73\times 10^{-9}$ \\
  & 1 & $-54.38401656$ & $1.57 \times 10^{-4}$ & $7.44 \times 10^{-4}$ & $4.77\times 10^{-6}$ & $7.39\times 10^{-9}$ \\
\hline
\multirow{2}{*}{3.0} 
  & 0 & $-50.18662388$ & $9.85 \times 10^{-5}$ & $1.07 \times 10^{-3}$ & $5.28\times 10^{-3}$ & $2.51\times 10^{-9}$ \\
  & 1 & $-47.80516396$ & $1.65 \times 10^{-4}$ & $1.23 \times 10^{-3}$ & $4.47\times 10^{-3}$ & $6.19 \times 10^{-9}$ \\
\hline
\multirow{2}{*}{2.0} 
  & 0 & $-35.90725762$ & $1.11 \times 10^{-4}$ & $1.10 \times 10^{-3}$ & $8.68 \times 10^{-3}$ & $2.42 \times 10^{-7}$ \\
  & 1 & $-35.52133640$ & $1.97 \times 10^{-4}$ & $1.01 \times 10^{-3}$ & $3.48\times 10^{-2}$  & $2.70 \times 10^{-7}$ \\
\hline
\multirow{2}{*}{1.0} 
  & 0 & $-26.86050464$ & $2.62 \times 10^{-5}$ & $1.38 \times 10^{-4}$ & $1.10\times 10^{-3}$ & $1.45 \times 10^{-6}$ \\
  & 1 & $-26.86046483$ & $2.46 \times 10^{-5}$ & $1.76 \times 10^{-4}$ &  & \\
\hline
\end{tabular}
\caption{Reference energies and relative eigenvalue errors of PEPS approximations for the first two eigenvalues ($\alpha=0,1$, corresponding to the ground state and first excited state) of the $4 \times 4$ TFI model \eqref{eq:TFI} at several values of $g$. We compare block-isoPEPS with $p=2$ at bond dimensions $\chi=4, \eta=8$ and $\chi=12, \eta=20$ against standard PEPS with bond dimensions $\chi=2$ and $\chi=4$. }
\label{tab:peps_isopeps_L4}
\end{table*}

\begin{figure}[t]
\centering
\includegraphics[scale=0.45]{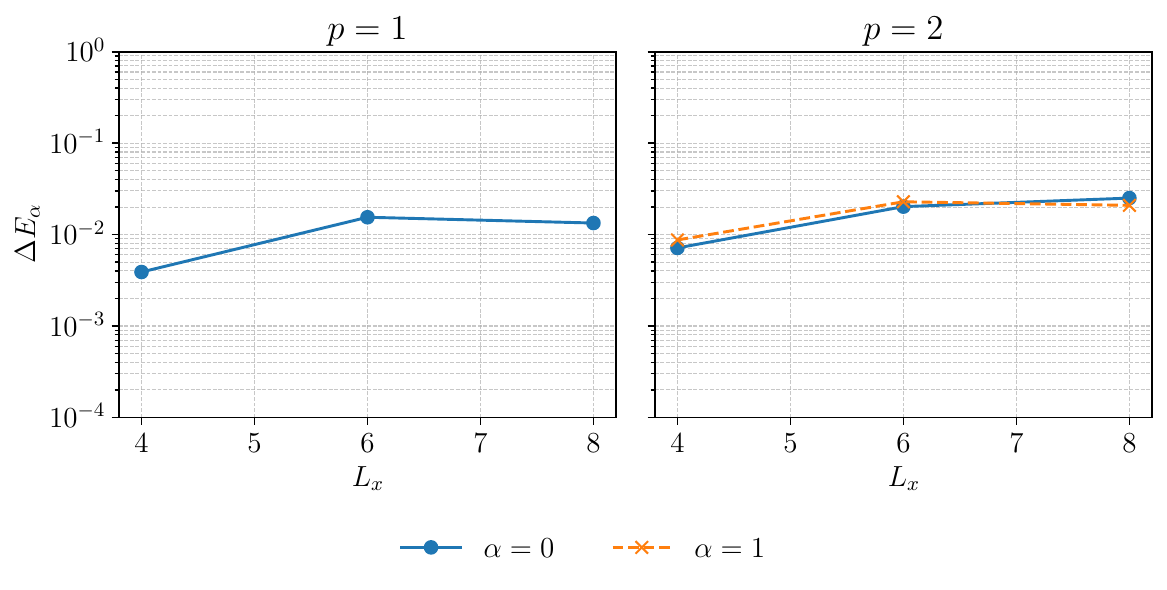}
\caption{
Relative eigenvalue error for the Heisenberg model \eqref{eq:heis} versus lattice side length $L_x$. Approximate eigenvalues were computed using $50$ iterations of block iso-PEPS subspace iteration with bond dimensions $\chi = 12$ and $\eta = 36$ and $\tau = 0.1$. Results are shown for block sizes $p = 1$ (left) and $p = 2$ (right).} 
\label{fig:Heis_XYZ_p2_versus_L}
\end{figure}

We also considered the antiferromagnetic Heisenberg model \eqref{eq:heis}. We computed (block) isoPEPS approximations of the ground state and first excited state using block sizes $p=1$ and $p=2$ with bond dimensions $\chi=12$ and $\eta=36$. We ran $50$ iterations of the isoPEPS subspace iteration algorithm with $\tau = 0.1$, and compared the results with reference energies obtained from DMRG. In \cref{fig:Heis_XYZ_p2_versus_L}, we show the relative error for lattices of side length $L_x=4,6,8$. Compared to the transverse-field Ising model, the ground and excited states of the Heisenberg model are significantly more challenging to represent using tensor networks, due to their higher entanglement. Increasing the block size in the isoPEPS representation from $p=1$ to $p=2$ does not significantly increase the error. 

\section{Conclusions} \label{sec:conclusions}

We introduced a block isometric tensor network method for computing excited states of two-dimensional quantum many-body systems via subspace iteration. Our approach extends the block (or state-averaged) MPS tensor network ansatz, widely used for excited-state calculations in one dimension, to the two-dimensional setting using projected entangled pair state tensor network with isometric constraints. Numerical experiments on the transverse field Ising and Heisenberg models demonstrate that the block-isoPEPS subspace iteration scales favorably with system size and can accurately capture low-energy excited states of spin systems across different parameter regimes. Furthermore, we demonstrated that the proposed method is capable of capturing excitations in systems where existing PEPS approaches based on the tangent space of the ground state are ineffective. 

Looking ahead, the block-isoPEPS method presented in this paper stands to benefit from several potential extensions. These include incorporating symmetries and employing DMRG-based algorithms specifically tailored to isometric PEPS \cite{Hyatt2020, Lin2022}, which may yield improved convergence. Another avenue of future research is extending block isometric tensor networks to three dimensions leveraging recent developments in 3D isometric tensor network algorithms \cite{Tepaske2021}. Finally, in future work we  plan to focus on parallel implementation and GPU acceleration to enhance computational efficiency and scalability of the proposed algorithms. 

\section*{Acknowledgments}
This material is based in part upon work supported by the U.S. Department of Energy, Office of Science, Office of Advanced Scientific Computing Research, Scientific Discovery through Advanced Computing (SciDAC) program through the FASTMath Institute and in part by the U.S. Department of Energy, Office of Science, Office of Advanced Scientific Computing Research's Applied Mathematics Competitive Portfolios program under U.S. Department of Energy Contract No. DE-AC02-05CH11231. This research used resources of the National Energy Research Scientific Computing Center, a DOE Office of Science User Facility supported by the Office of Science of the U.S. Department of Energy under Contract No. DE-AC02-05CH11231 using NERSC award ASCR-ERCAPm1027.
R.C. was supported by the US Department of Energy, Office of Science, Office of Advanced Scientific Computing Research, under Award Number DE-SC0025572. 
The authors thank Garnet Kin-Lic Chan for useful discussions.

\bibliographystyle{elsarticle-num}
\bibliography{refs}

\appendix

\section{Finite PEPS tangent-space excitation method} \label{sec:peps_tangent_space}
The tangent-space formalism was originally introduced for matrix product states (MPS) by \"Ostlund and Rommer~\cite{Ostlund1995}, and subsequently generalized to infinite MPS by Haegeman et al.~\cite{Haegeman2012}. Its extension to infinite PEPS was later developed by Vanderstraeten et al.~\cite{Vanderstraeten2015}.  
The tangent-space method is most commonly employed in the context of translation-invariant system, where the excitation ansatz is explicitly constructed as an eigenstate of the translation operator. In order to facilitate a meaningful comparison with finite-size isoPEPS calculations, we adapt this framework to the finite-system setting.

We begin with a fully optimized ground state represented by a PEPS \eqref{eq:peps}. The tangent space of this ground state is spanned by tensors $T_{\boldsymbol r}(B_{\boldsymbol r})$, defined by replacing a single tensor at lattice position $\boldsymbol{r}=(i,j)$ in the ground-state PEPS $T$ with another tensor $B_{\boldsymbol r}$. For convenience, we relabel these basis vectors as $T(B_m)$, where the index $m$ encodes the site coordinate. This set of tensors is overcomplete and mutually non-orthogonal. The corresponding Gram matrix is 
\begin{equation}
N_{mn} = \left\langle T(B_{m}), T(B_{n}) \right\rangle,
\end{equation}
and the projected Hamiltonian matrix is given by 
\begin{equation}
H_{mn} = \left\langle T(B_{m}), H T(B_{n})\right\rangle,
\end{equation}
where $\langle \bullet, \bullet \rangle$ denotes the Frobenius inner product of tensors obtained by contracting over all shared physical indices. 

Since $T(B_n)$ is a linear function of $B_n$, both $N_{mn}$ and $H_{mn}$ take the form of quadratic expressions in $B_n$ 
\begin{align}
N_{mn} &= B_m^{\ast} \mathcal{N} B_n, \notag \\ 
H_{mn} &= B_m^{\ast} \mathcal{H} B_n,
\end{align}
where $\mathcal{N}$ and $\mathcal{H}$ are independent of the particular choice of basis tensor $B_n$. This structure implies that the elements of the effective Gram matrix $N_{mn}$ can be computed by first deriving the vector form $\mathcal{N} B_n$ and then contracting it with $B_m^{\ast}$ as $
N_{mn} = B_m^{\ast} (\mathcal{N} B_n)$, with an analogous expression holding for $H_{mn}$. The corresponding matrix–vector products can be expressed as 
\begin{equation}
\begin{aligned}
\mathcal{N} B_n &= \frac{\partial}{\partial B_n^{\ast}} \big( B_n^{\ast} \mathcal{N} B_n \big)  \\ 
\mathcal{H} B_n &= \frac{\partial}{\partial B_n^{\ast}} \big( B_n^{\ast} \mathcal{H} B_n \big).
\end{aligned}
\end{equation}
These derivatives can be efficiently evaluated using automatic differentiation, as discussed in Ref.~\cite{Ponsioen2022}. Once the effective Hamiltonian $H_{mn}$ and norm matrix $N_{mn}$ are obtained, the excitation spectrum is determined by solving the generalized eigenvalue problem
\begin{equation}
\sum_n H_{mn} v_{np} = E_p \sum_n N_{mn} v_{np}.
\end{equation}
It is important to note that, because the basis is overcomplete due to the gauge invariance of PEPS, $N_{mn}$ contains zero modes. These must be projected out prior to solving the eigenvalue problem, which may otherwise lead to numerical instabilities. In contrast, in the isoPEPS construction the block states are orthogonal by design, and this issue does not arise.

The main computational cost of the tangent-space algorithm lies in evaluating the elements of $H_{mn}$, each of which requires a tensor contraction scaling as $\bigO(L_xL_y \chi^{10})$ for every basis state. Since the number of basis states scales as $O(L_xL_y \chi^4)$, the overall complexity is $\bigO(L_x^2L_y^2 \chi^{14})$. This cost can be reduced by restricting to a finite subset $M$ of basis states in the tangent space, yielding a scaling of $\bigO(M L_xL_y \chi^{10})$, although it still remains higher. 

\end{document}